\documentclass[12pt,english]{smfart}

\usepackage{smfthm}
\usepackage{smfenum}
\usepackage{amssymb}
\usepackage{amscd}

\newcommand{\Qp}{\mathbf{Q}_p}
\newcommand{\Zp}{\mathbf{Z}_p}
\newcommand{\Cp}{\mathbf{C}_p}

\newcommand{\eps}{\varepsilon}
\newcommand{\epsn}{\eps^{(n)}}
\newcommand{\ra}{\rightarrow}

\newcommand{\Qpbar}{\overline{\mathbf{Q}}_p}

\newcommand{\on}{\operatorname}
\newcommand{\OO}{\mathcal{O}}

\renewcommand{\hat}{\widehat}

\renewcommand{\phi}{\varphi}
\renewcommand{\projlim}{\varprojlim}
\newcommand{\ZZ}{\mathbf{Z}}

\newcommand{\btrig}[2]{\widetilde{\mathbf{B}}^{\dagger #1}_{\mathrm{rig} #2}}

\newcommand{\bnrig}[2]{\mathbf{B}^{\dagger #1}_{\mathrm{rig} #2}}

\newcommand{\btrigplus}[1]{\widetilde{\mathbf{B}}^{+}_{\mathrm{rig} #1}}

\newcommand{\dnrig}[1]{\mathbf{D}^{\dagger #1}_{\mathrm{rig}}}
\newcommand{\dhol}{\mathbf{D}^+_{\mathrm{rig}}}

\newcommand{\bmax}{\mathbf{B}_{\mathrm{max}}}
\newcommand{\bcris}{\mathbf{B}_{\mathrm{cris}}}

\newcommand{\bdR}{\mathbf{B}_{\mathrm{dR}}}
\newcommand{\btdag}[1]{\widetilde{\mathbf{B}}^{\dagger #1}}
\newcommand{\bdag}[1]{\mathbf{B}^{\dagger #1}}
\newcommand{\bplus}{\mathbf{B}^+}
\newcommand{\btplus}{\widetilde{\mathbf{B}}^+}
\newcommand{\bt}{\widetilde{\mathbf{B}}}

\newcommand{\aplus}{\mathbf{A}^+}
\newcommand{\atplus}{\widetilde{\mathbf{A}}^+}
\newcommand{\at}{\widetilde{\mathbf{A}}}
\newcommand{\bhol}[1]{\mathbf{B}^+_{\mathrm{rig} #1}}

\newcommand{\e}{\mathbf{E}}

\newcommand{\etplus}{\widetilde{\mathbf{E}}^+}
\newcommand{\et}{\widetilde{\mathbf{E}}}

\newcommand{\dcris}{\mathbf{D}_{\mathrm{cris}}}
\newcommand{\ddR}{\mathbf{D}_{\mathrm{dR}}}
\newcommand{\dfont}{\mathbf{D}}

\renewcommand{\ddag}[1]{\mathbf{D}^{\dagger #1}}
\newcommand{\dsen}{\mathbf{D}_{\mathrm{Sen}}}

\newcommand{\vale}{v_\mathbf{E}}

\title{Bloch and Kato's exponential map: three explicit formulas}

\author{Laurent Berger}
\address{Harvard Dept of Mathematics  \\
      One Oxford Street \\
      Cambridge, MA 02138-2901 \\ USA}
\email{laurent@math.harvard.edu}
\urladdr{www.math.harvard.edu/\~{}laurent}

\thanks{This research was partially conducted by the author for the
  Clay Mathematical Institute}


\date{\today}
\begin{document}

\begin{abstract}
The purpose of this article is to give formulas for Bloch-Kato's
exponential map and its dual for an absolutely crystalline 
$p$-adic representation $V$, in terms of the $(\phi,\Gamma)$-module
associated to that representation. 
As a corollary of these computations, we can give a very simple
(and slightly improved)
description of Perrin-Riou's exponential map (which interpolates
Bloch-Kato's exponentials for $V(k)$). This
new description directly implies Perrin-Riou's
reciprocity formula.
\end{abstract}

\maketitle
\setcounter{tocdepth}{2}
\tableofcontents

\setlength{\baselineskip}{18pt}

\section*{Introduction}
Let $p$ be a prime number, and let $V$ be a $p$-adic representation of
$\on{Gal}(\overline{K}/K)$ where $K$ is a finite extension of
$\Qp$. Such objects arise (for example) as the \'etale cohomology of
algebraic varieties, hence their interest in arithmetic algebraic
geometry. 

Let $\bcris$ and $\bdR$ be the rings of periods of Fontaine, and
let $\dcris(V)$ and $\ddR(V)$ be the invariants attached to $V$ by
Fontaine's construction. Bloch and Kato have defined (in 
\cite{BK91}),
for a de Rham representation $V$, an
``exponential'' map, $\exp_{K,V} : \ddR(V) / \on{Fil}^0 \ddR(V) \ra H^1(K,V)$. 
It is obtained by tensoring the so-called fundamental exact sequence:
\[ 0 \ra \Qp \ra \bcris^{\phi=1} \ra \bdR/\bdR^+ \ra 0 \]
with $V$ and taking the invariants under the action of $G_K$. The
exponential map is then the connecting homomorphism 
$\ddR(V) / \on{Fil}^0 \ddR(V) \ra H^1(K,V)$.

The reason
for their terminology is the following (cf \cite[3.10.1]{BK91}):
if $G$ is a formal Lie group of finite
height over $\OO_K$, and $V = \Qp \otimes_{\Zp} T$ where
$T$ is the $p$-adic Tate module of $G$, then $V$ is a de Rham
representation and $\ddR(V)/\on{Fil}^0 \ddR(V)$ is identified with the
tangent space $\tan(G(K))$ of $G(K)$. 
In this case, we have a commutative diagram:
\[ \begin{CD}
\tan(G(K)) @>{\exp_G}>> \mathbf{Q} \otimes_{\ZZ} G(\OO_K) \\
@VVV @V{\delta_G}VV \\
\ddR(V)/\on{Fil}^0 \ddR(V) @>{\exp_{K,V}}>> H^1(K,V), \end{CD} \]
where $\delta_G$ is the Kummer map, the upper $\exp_G$ is the usual
exponential map, and the lower $\exp_{K,V}$ is Bloch-Kato's exponential map.

The cup product $\cup : H^1(K,V) \times H^1(K,V^*(1)) \ra H^2(K,\Qp(1))
\simeq \Qp$ defines a perfect pairing, which we can use 
(by dualizing twice) to define
Bloch and Kato's dual exponential map $\exp^*_{K,V^*(1)} : H^1(K,V)
\ra \on{Fil}^0 \ddR(V)$.

Finally, Perrin-Riou has constructed in \cite{BP94}
a period map $\Omega_{V,h}$ which
interpolates the $\exp_{K,V(k)}$ as $k$ runs over the positive
integers. It is a crucial ingredient in the construction of $p$-adic
$L$ functions, and is a vast generalization of Coleman's map.

The goal of this article is to give formulas for $\exp_{K,V}$, 
$\exp^*_{K,V^*(1)}$, and $\Omega_{V,h}$ in terms of the
$(\phi,\Gamma)$-module associated to $V$ by Fontaine. As a corollary,
we recover a theorem of Colmez which states that Perrin-Riou's map 
interpolates the $\exp^*_{K,V^*(1-k)}$ as $k$ runs over the negative
integers. This is equivalent to Perrin-Riou's 
conjectured reciprocity formula.
Our construction of $\Omega_{V,h}$ is actually a slight improvement
over Perrin-Riou's (one does not have to kill the $\Lambda$-torsion, 
see remark \ref{no_kill}).

We refer the reader to the text itself for a statement of the actual
formulas (theorems \ref{exp_bk}, \ref{exp_duale} and \ref{bpr_const}) 
which are rather technical.

This article does not really 
contain any new results, and it is mostly a
re-interpretation of formulas of Cherbonnier-Colmez
(for the dual exponential map), 
and of Colmez and Benois 
(for Perrin-Riou's map)
in the language of the author's article ``$p$-adic
representations and differential equations''.  

\renewcommand{\thesection}{\Roman{section}}
\section{Periods of $p$-adic representations}

\Subsection{Notations}\label{notations}
Throughout this article, $k$ will denote a finite field of
characteristic $p > 0$, so that $F=W(k)[1/p]$ is a finite unramified
extension of $\Qp$. Let $\Qpbar$ be the algebraic closure of $\Qp$, 
let $K$ be a finite totally ramified extension of $F$, and
let $G_K = \on{Gal}(\Qpbar/K)$ be the absolute Galois group of $K$. 
Let $\mu_{p^n}$ be the group of $p^n$th 
roots of unity; for every $n$, we will choose a generator 
$\epsn$ of $\mu_{p^n}$, with the additional requirement that
$(\epsn)^p=\eps^{(n-1)}$. This makes $\projlim_n \epsn$ into a generator
of $\projlim_n \mu_{p^n} \simeq \Zp(1)$. We set $K_n=K(\mu_{p^n})$ and
$K_{\infty}=\cup_{n=0}^{+\infty} K_n$. Recall that the cyclotomic
character $\chi: G_K \ra \Zp^*$ is defined by the relation:
$g(\epsn)=(\epsn)^{\chi(g)}$ for all $g \in G_K$. The kernel of the
cyclotomic character is $H_K = \on{Gal}(\Qpbar/K_{\infty})$, and
$\chi$ therefore identifies $\Gamma_K = G_K / H_K$ with
an open subgroup of $\Zp^*$. The
torsion subgroup of $\Gamma_K$ will be denoted by $\Delta_K$.
We also set $\Gamma_n=\on{Gal}(K_{\infty}/K_n)$. When $p \neq 2$ and
$n \geq 1$ (or $p=2$ and $n \geq 2$), $\Gamma_n$ is torsion free. 
If $x \in 1+p\Zp$, then there exists $k \geq 1$ such that 
$\log_p(x) \in p^k \Zp^*$ and we'll write $\log_p^0(x) =
\log_p(x)/p^k$. 

The completed group algebra of $\Gamma_K$ is
$\Lambda=\Zp[[\Gamma_K]] \simeq \Zp[\Delta_K] \otimes_{\Zp} 
\Zp[[\Gamma_1]]$, and we set $\mathcal{H}(\Gamma_K)=
\Qp[\Delta_K] \otimes_{\Qp} \mathcal{H}(\Gamma_1)$ where 
$\mathcal{H}(\Gamma_1)$ is the set of $f(\gamma-1)$ with $\gamma \in
\Gamma_1$ and where $f(X) \in \Qp[[X]]$ is convergent on the $p$-adic open
unit disk. Examples of elements of $\mathcal{H}(\Gamma_K)$ are the
$\nabla_i$ (which are Perrin-Riou's $\ell_i$'s), defined by $\nabla_i
= \log(\gamma)/\log_p(\chi(\gamma))-i$. We will also use the
operator $\nabla_0/(\gamma_n-1)$, where $\gamma_n$ is a topological
generator of $\Gamma_n$ (see \cite[4.1]{LB02}). 
It is defined by the formula:
\[ \frac{\nabla_0}{\gamma_n-1} = 
\frac{\log(\gamma_n)}{\log_p(\chi(\gamma_n))(\gamma_n-1)} =
\frac{1}{\log_p(\chi(\gamma_n))} \sum_{n \geq 1}
\frac{(1-\gamma_n)^{n-1}}{n}, \]
or equivalently by 
\[ \frac{\nabla_0}{\gamma_n-1} = \lim_{\Gamma \ni \eta \ra 1}
\frac{\eta-1}{\gamma_n-1}
\frac{1}{\log_p(\chi(\eta))}.\]
It is easy to see that 
$\nabla_0/(\gamma_n-1)$ acts on $F_n$ by $1/\log_p(\chi(\gamma_n))$.

The algebra acts $\mathcal{H}(\Gamma_K)$ on
$\bhol{,F}$ and one can check that 
\[ \nabla_i  = t\frac{d}{dt} -i = \log(1+\pi)\partial-i, 
\quad\text{where}\quad \partial = (1+\pi)\frac{d}{d\pi}. \]
In particular, $\nabla_0 \bhol{,F} \subset t \bhol{,F}$ and if $i \geq
1$, then 
$\nabla_{i-1} \circ \cdots \circ \nabla_0 \bhol{,F} \subset t^i \bhol{,F}$.

A $p$-adic representation $V$ is a finite dimensional $\Qp$-vector
space with a continuous linear action of $G_K$. It is easy to see that
there is always a $\Zp$-lattice of $V$ which is stable by the
action of $G_K$, and such lattices will be denoted by $T$.
The main strategy (due to Fontaine, see for example \cite{Bu88sst})
for studying $p$-adic representations of a group $G$ 
is to construct topological
$\Qp$-algebras $B$, endowed with an action of $G$ and 
some additional structures so that if
$V$ is a $p$-adic representation, then 
$D_B(V)=(B \otimes_{\Qp} V)^G$ is a $B^G$-module which
inherits these structures, and so that the 
functor $V \mapsto D_B(V)$ gives interesting
invariants of $V$. We say that a $p$-adic 
representation $V$ of $G$ is $B$-admissible if we
have $B\otimes_{\Qp} V \simeq B^d$ as $B[G]$-modules.

\Subsection{$p$-adic Hodge theory}
In this paragraph, we will recall the definitions of Fontaine's rings
of periods. One can find some of these constructions in \cite{Bu88per}.
Let \[ \et=\projlim_{x\mapsto x^p} \Cp 
=\{ (x^{(0)},x^{(1)},\cdots) \mid (x^{(i+1)})^p = x^{(i)} \}, \]
and let $\etplus$ be the set of $x \in \et$ 
such that $x^{(0)} \in  \OO_{\Cp}$. 
If $x=(x^{(i)})$ and $y=(y^{(i)})$ are two elements of $\et$,
we define their sum $x+y$ and their product $xy$ by:
\[ (x+y)^{(i)}= \lim_{j \ra \infty} (x^{(i+j)}+y^{(i+j)})^{p^j} 
\quad\text{and}\quad
(xy)^{(i)}=x^{(i)}y^{(i)}, \] which 
makes $\et$ into an algebraically closed field of
characteristic $p$. If $x=(x^{(n)})_{n \geq 0} 
\in \et$, let $\vale(x)=v_p(x^{(0)})$.
This is a valuation on $\et$ for which $\et$ is 
complete; the ring of integers of $\et$ is
$\etplus$. Let $\atplus$ be the ring $W(\etplus)$ 
of Witt vectors with coefficients in
$\etplus$ and let $\btplus=\atplus[1/p]=\{ 
\sum_{k\gg -\infty} p^k [x_k],\ x_k \in \etplus
\}$ where $[x] \in \atplus$ is the Teichm{\"u}ller 
lift of $x \in \etplus$.  This ring is
endowed with a map $\theta: \btplus \ra \Cp$ defined 
by the formula 
$\theta(\sum_{k\gg -\infty} p^k[x_k] )
=\sum_{k\gg -\infty} p^k x_k^{(0)}$.
Let $\eps=(\eps^{(i)})\in\etplus$ where $\eps^{(n)}$ 
is defined above, and define $\pi=[\eps]-1$, $\pi_1=
[\eps^{1/p}]-1$, $\omega=\pi/\pi_1$
and $q=\phi(\omega)=\phi(\pi)/\pi$. 
One can show that $\ker(\theta : \atplus \ra \atplus)$
is the principal ideal generated by $\omega$.

Notice that $\eps \in \etplus$ and that 
$\vale(\eps-1)=p/(p-1)$. We define $\e_F=k((\eps-1))$ and
$\e$ as the separable closure of $\e_F$ in $\et$, 
as well as $\e^+=\e \cap \etplus$, the ring
of integers of $\e$. By definition, $\e$ is 
separably closed and $\et$ is the completion of
its purely inseparable closure.

The ring $\bdR^+$ is defined to be the 
completion of $\btplus$
for the $\ker(\theta)$-adic topology: \[
\bdR^+=\projlim_{n\geq 0} \btplus/(\ker(\theta)^n). \] 
It is a discrete valuation ring, whose maximal ideal
is $\ker(\theta)$; the series which 
defines $\log([\eps])$ converges in $\bdR^+$ to an element $t$, which
is a generator of the maximal ideal, so 
that $\bdR=\bdR^+[1/t]$ is a field, endowed with an action of $G_F$
and a filtration defined by $\on{Fil}^i(\bdR)=t^i \bdR^+$ for $i \in \ZZ$.

We say that a representation $V$ of $G_K$ 
is de Rham if it is $\bdR$-admissible which is
equivalent to the fact that the filtered $K$-vector 
space $\ddR(V)=(\bdR\otimes_{\Qp} V)^{G_K}$
is of dimension $d=\dim_{\Qp}(V)$. 

The ring $\bmax^+$ is defined as being
\[ \bmax^+= \{ \sum_{n \geq 0} a_n \frac{\omega^n}{p^n}
\text{ where $a_n\in \btplus$ is sequence converging to $0$} \}, \]
and $\bmax=\bmax^+[1/t]$. One could replace 
$\omega$ by any generator of $\ker(\theta)$. The
ring $\bmax$ injects canonically into $\bdR$ and, 
in particular, it is endowed with the induced
Galois action and filtration, as well as 
with a continuous Frobenius $\phi$, extending the
map $\phi: \atplus \ra \atplus$ coming 
from $x \mapsto x^p$ in $\etplus$. Let us point out that
$\phi$ does not extend continuously 
to $\bdR$. One also sets $\btrigplus{}=\cap_{n=0}^{+\infty}
\phi^n(\bmax^+)$.

We say that a representation $V$ 
of $G_K$ is crystalline if it is $\bmax$-admissible or
(which is the same) $\btrigplus{}[1/t]$-admissible 
(the periods of crystalline representations
live in finite dimensional $F$-vector 
subspaces of $\bmax$, stable by $\phi$, and so in fact
in $\cap_{n=0}^{+\infty}
\phi^n(\bmax^+)[1/t]$); this is equivalent 
to requiring that the $F$-vector space
\[ \dcris(V)=(\bmax \otimes_{\Qp} V)^{G_K}=(\btrigplus{}[1/t]
\otimes_{\Qp} V)^{G_K} \] 
be of dimension $d=\dim_{\Qp}(V)$. 
Then $\dcris(V)$ is endowed with a Frobenius $\phi$ 
induced by that of $\bmax$,  
and $(\bdR\otimes_{\Qp} V)^{G_K}=\ddR(V)=
K \otimes_F \dcris(V)$ so that a crystalline 
representation is also de Rham and $K \otimes_F \dcris(V)$ is a
filtered $K$-vector space.

If $V$ is a $p$-adic representation, we say 
that $V$ is Hodge-Tate, with Hodge-Tate weights
$h_1, \cdots, h_d$, if we have a decomposition $\Cp\otimes_{\Qp} V \simeq
\oplus_{j=1}^d \Cp(h_j)$. In this case, we see that 
$(\Cp\otimes_{\Qp} V)^{H_K} \simeq 
\oplus_{j=1}^d \hat{K}_{\infty} (h_j)$ and one can show
that $\dsen(V)=(\Cp\otimes_{\Qp} V)^{H_K}_{\mathrm{fin}}$, 
which is by definition the union of the finite
dimensional $K_{\infty}$-vector subspaces  
of $(\Cp\otimes_{\Qp} V)^{H_K}$ which are stable by $\Gamma_K$,
is equal to $\oplus_{j=1}^d K_{\infty} (h_j)$.

The $K_{\infty}$-vector space $\dsen(V)$ is
endowed with a residual action of $\Gamma_K$ 
and if $\gamma \in \Gamma_K$ is sufficiently
close to $1$, then the series of operators 
$\log(\gamma)/\log_p(\chi(\gamma))$ converges to a
$K_{\infty}$-linear operator $\nabla_V : 
\dsen(V) \ra \dsen(V)$ which does not depend on the
choice of $\gamma$, and which is 
diagonalizable with eigenvalues $h_1, \cdots, h_d$. We will
say that $V$ is positive if its Hodge-Tate 
weights are negative (the definition of the sign
of the Hodge-Tate weights is unfortunate; some people change the
sign and talk about geometrical weights). 
By using the map $\theta: \bdR^+ \ra \Cp$,
it is easy to see that the
Hodge-Tate weights of $V$ are those integers $h$ such that
$\on{Fil}^{-h} \ddR(V) \neq \on{Fil}^{-h+1} \ddR(V)$.

\Subsection{$(\phi,\Gamma)$-modules}
See \cite{Fo91} for a reference. 
Let $\at$ be the ring of Witt vectors 
with coefficients in $\et$ and $\bt=\at[1/p]$. Let
$\mathbf{A}_F$ be the completion of 
$\OO_F[\pi,\pi^{-1}]$ in $\at$ for this ring's
topology, which is also the completion 
of $\OO_F[[\pi]][\pi^{-1}]$ for the $p$-adic topology
($\pi$ being small in $\at$).
This is a discrete valuation ring whose 
residual field is $\e_F$. Let $\mathbf{B}$ be the
completion for the $p$-adic topology of 
the maximal unramified extension of
$\mathbf{B}_F=\mathbf{A}_F[1/p]$ in $\bt$. 
We then define $\mathbf{A}=\mathbf{B}
\cap \at$
and $\aplus=\mathbf{A} \cap \atplus$.
These rings are endowed with an action 
of Galois and a Frobenius deduced from those on $\et$.
We set $\mathbf{A}_K=\mathbf{A}^{H_K}$ and
$\mathbf{B}_K=\mathbf{A}_K[1/p]$. When $K=F$, the two definitions are the same.
One also sets $\bplus=\aplus[1/p]$, and $\bplus_F=(\bplus)^{H_F}$ as well as
$\aplus_F=(\aplus)^{H_F}$.

If $V$ is a $p$-adic representation 
of $G_K$, let $\dfont(V)=(\mathbf{B} \otimes_{\Qp}
V)^{H_K}$. We know by \cite{Fo91} that 
$\dfont(V)$ is a $d$-dimensional $\mathbf{B}_K$-vector
space with a Frobenius and a residual action 
of $\Gamma_K$ which commute (it is a
$(\phi,\Gamma_K)$-module) and that one can recover $V$ by 
the formula
$V=(\dfont(V) \otimes_{\mathbf{B}_K} \mathbf{B})^{\phi=1}$. 

The field $\mathbf{B}$ is a totally ramified extension (because the
residual extension is purely inseparable) of degree $p$ of
$\phi(\mathbf{B})$.
The Frobenius map $\phi: \mathbf{B} \ra \mathbf{B}$ 
is injective but therefore not surjective, but we can define a left
inverse for $\phi$, which will play a major role in the sequel. We set:
$\psi(x)=\phi^{-1}(p^{-1} \on{Tr}_{\mathbf{B}/\phi(\mathbf{B})}(x))$.

Let us now set $K=F$. We say that a $p$-adic 
representation $V$ of $G_F$ is of finite height 
if $\dfont(V)$ has a basis over $\mathbf{B}_F$ 
made up of elements of $\dfont^+(V) = (\bplus
\otimes_{\Qp} V)^{H_F}$. A result of Fontaine 
(\cite{Fo91} or \cite[III.2]{Co99})
shows that $V$ is of finite height if and only 
if $\dfont(V)$ has a sub-$\bplus_F$-module
which is free of finite rank $d$, and stable 
by $\phi$. Let us recall the main result (due to
Colmez, see \cite{Co99} or also \cite[theorem 3.10]{LB02}) 
regarding crystalline representations of $G_F$:
if $V$ is a crystalline representation of $G_F$, 
then $V$ is of finite height. 

If $K \neq F$ or if $V$ is no longer crystalline, then it is no longer
true in general that $V$ is of finite height, but it is still
possible to say something about the periods of $V$. Every element $x
\in \bt$ can be written in a unique way as $x = \sum_{k \gg -\infty}
p^k [x_k]$, where $x_k \in \et$. For $r > 0$, let us set:
\[\btdag{,r}=\left\{ x \in \bt,\ 
\lim_{k \ra +\infty} \vale(x_k)+\frac{pr}{p-1}k 
= +\infty \right\}. \]
This makes $\btdag{,r}$ into an intermediate ring between $\btplus$
and $\bt$. Let us set 
$\bdag{,r} = \mathbf{B} \cap \btdag{,r}$,
$\btdag{} =\cup_{r \geq 0} \btdag{,r}$, and
$\bdag{} = \cup_{r\geq 0} \bdag{,r}$. 
If $R$ is any of the above rings, then by definition $R_K=R^{H_K}$.

We say that a $p$-adic representation $V$ is overconvergent if
$\dfont(V)$ has a basis over $\mathbf{B}_K$ made up of elements of
$\ddag{}(V) = (\bdag{} \otimes_{\Qp} V)^{H_K}$. The main result on
the overconvergence of
$p$-adic representations of $G_K$ is the following (cf  \cite{CC98,CC99}):
\begin{theo}
Every $p$-adic representation $V$ of $G_K$ is overconvergent, that is
there exists $r(V)$ such that 
$\dfont(V)=\mathbf{B}_K \otimes_{\bdag{,r(V)}_K} \ddag{,r(V)}(V)$.
\end{theo}

The terminology ``overconvergent'' can be explained by the following
proposition, which describes the rings $\bdag{,r}_K$. Let
$e_K=[K_\infty: F_\infty]$:

\begin{prop}
Let $\mathcal{B}_F^\alpha$ be the set of power series 
$f(X) = \sum_{k \in \ZZ} a_k
X^k$ such that $a_k$ is a bounded sequence of elements of $F$, and
such that $f(X)$ is holomorphic on the $p$-adic annulus 
$\{ p^{-1/\alpha} \leq |T| < 1 \}$.

There exists $r(K)$ and $\pi_K \in \bdag{,r(K)}_K$ such that if $r
\geq r(K)$, then the map $f \mapsto f(\pi_K)$ from 
$\mathcal{B}_F^{e_K r}$ to $\bdag{,r}_K$ is an isomorphism.
If $K=F$, then one can take $\pi_F=\pi$.
\end{prop}

\Subsection{$p$-adic representations and differential equations}
We shall now recall some of the results of \cite{LB02}, which allow us
to recover $\dcris(V)$ from the $(\phi,\Gamma)$-module of $V$. Let 
$\mathcal{H}_F^\alpha$ be the set of power series 
$f(X) = \sum_{k \in \ZZ} a_k
X^k$ such that $a_k$ is a sequence 
(not necessarily bounded)
of elements of $F$, and
such that $f(X)$ is holomorphic on 
the $p$-adic annulus $\{ p^{-1/\alpha} \leq |T| < 1 \}$.

For $r \geq r(K)$, define $\bnrig{,r}{,K}$ as the set of $f(\pi_K)$
where $f(X) \in \mathcal{H}_F^{e_K r}$. Obviously, 
$\bdag{,r}_K \subset \bnrig{,r}{,K}$. If $V$ is a $p$-adic
representation, let $\dnrig{,r}(V) = \bnrig{,r}{,K} \otimes_{\bdag{,r}_K}
\ddag{,r}(V)$. 

One of the main technical tools 
of \cite{LB02} is the construction of a large
ring $\btrig{}{}$, which contains $\btrigplus{}$ and $\btdag{}$, so that 
$\dcris(V) \subset (\btrig{}{}[1/t] \otimes_{\Qp} V)^{G_K}$ and 
$\dnrig{}(V)[1/t] \subset 
(\btrig{}{}[1/t] \otimes_{\Qp} V)^{H_K}$. We then have
(cf \cite[theorem 3.6]{LB02}): $\dcris(V) = (\dnrig{}(V)[1/t])^{\Gamma_F}$.

Let us now return to the case when $K=F$ and $V$ is a crystalline
representation of $G_F$. In this case, one can give a more precise
result. Let $\bhol{,F}$ be the set of $f(\pi)$ where $f(X)=\sum_{k
  \geq 0} a_k X^k$ where $a_k \in F$, and $f(X)$ is holomorphic on the
$p$-adic open unit disk. 
Set $\dhol(V) = \bhol{,F} \otimes_{\bplus_F} \dfont^+(V)$.
One can then show (cf \cite{LBcr}):
$\dcris(V) = (\dhol(V)[1/t])^{\Gamma_F}$.

\Subsection{Construction of cocycles}
The purpose of this paragraph is to recall the constructions of
\cite[I.5]{CC99} and extend them a little bit. 
In this paragraph and in the next, $V$ will be an arbitrary $p$-adic
representation of $G_K$.
Recall that in
loc. cit., a map $h^1_{K,V} : \dfont(V)^{\psi=1} \ra H^1(K,V)$ was
constructed, and that (when $\Gamma_K$ is torsion free at least)
it gives rise to an exact sequence:
\[ \begin{CD}
0 @>>> \dfont(V)^{\psi=1}_{\Gamma_K} @>{h^1_{K,V}}>> H^1(K,V) @>>> \left( 
\frac{\dfont(V)}{\psi-1} \right)^{\Gamma_K} @>>> 0. \end{CD} \]
We shall extend $h^1_{K,V}$ to a map $h^1_{K,V}:
\dnrig{}(V)^{\psi=1} \ra H^1(K,V)$. 
We will first need a few facts about the ring of periods $\btrig{}{}$
and the modules $\dnrig{,r}(V)$.

\begin{lemm}\label{inv_gam}
If $\gamma \in \Gamma_K$, then $1-\gamma: \dnrig{,r}(V) \ra
\dnrig{,r}(V)$ is an isomorphism.
\end{lemm}

\begin{proof}
Recall that $\bnrig{,r}{,K}$ is the completion of $\bdag{,r}_K$ for
the Fr\'echet topology (see \cite[2.6]{LB02}), so 
that $\dnrig{,r}(V)$ is the completion of $\ddag{,r}(V)$ for
the Fr\'echet topology.
The lemma then follows from the fact that
by \cite[II.6.1]{CC98},
if $\gamma \in \Gamma_K$, then $1-\gamma: \ddag{,r}(V) \ra
\ddag{,r}(V)$ is a bi-continuous isomorphism for the Fr\'echet topology.
\end{proof}

\begin{lemm}\label{inv_phi}
There is an exact sequence:
\[ \begin{CD} 
0 @>>> \Qp @>>> \btrig{}{} @>{1-\phi}>> \btrig{}{} @>>> 0 
\end{CD} \] 
\end{lemm}

\begin{proof}
We'll start with the easiest part, namely the fact that
$(\btrig{}{})^{\phi=1}=\Qp$. If $x \in (\btrig{}{})^{\phi=1}$, then
\cite[prop 3.2]{LB02} shows that actually $x \in (\btrigplus{})^{\phi=1}$,
and the latter space is well-known to be $=\Qp$.

Let us now show that if $\alpha \in \btrig{}{}$, then there exists
$\beta \in \btrig{}{}$ such that $(1-\phi)\beta=\alpha$. By
\cite[lemma 2.18]{LB02}, one can write $\alpha=\alpha^++\alpha^-$ with
$\alpha^+ \in \btrigplus{}$ and $\alpha^- \in \btdag{}{}$. It is
therefore enough to show that $1-\phi: \btrigplus{} \ra \btrigplus{}$
and $1-\phi: \btdag{}{} \ra \btdag{}{}$ are surjective. 

The first assertion follows from the 
facts that $1-\phi : \bcris^+ \ra \bcris^+$ 
is surjective (see \cite[thm 5.3.7, ii]{Bu88per}) and that
$\btrigplus{} = \cap_{n=0}^{+\infty} \phi^n(\bcris^+)$.

The second assertion follows from the facts that $1-\phi: \bt \ra \bt$
is surjective and that if $\beta \in \bt$ is such that $(1-\phi)\beta
\in \btdag{}{}$, then $\beta \in \btdag{}{}$.
\end{proof}

If $K$ and $n$
are such that $\Gamma_n$ is torsion-free, then
we will construct maps $h^1_{K_n,V}$ such that 
$\on{cor}_{K_{n+1}/K_n} \circ h^1_{K_{n+1},V} = h^1_{K_n,V}$. If 
$\Gamma_n$ is no longer torsion free, we'll therefore
define $h^1_{K_n,V}$ by the
relation $h^1_{K_n,V} = \on{cor}_{K_{n+1}/K_n} \circ h^1_{K_{n+1},V}$.
In the following proposition, we therefore assume that $\Gamma_K$ is
torsion free (and therefore procyclic), 
and we let $\gamma$ denote a topological generator of $\Gamma_K$.
If $M$ is a $\Gamma_K$-module, it is customary to write 
$M_{\Gamma_K}$ for $M/ \on{im}(\gamma-1)$. 

\begin{prop}\label{cocyle}
If $y \in \dnrig{}(V)^{\psi=1}$, and $b \in \btrig{}{} \otimes_{\Qp} V$ 
is a solution of the equation $(\gamma-1)(\phi-1)b=(\phi-1)y$, then
the formula \[ h^1_{K,V}(y) = 
\log_p^0(\chi(\gamma)) \left[ \sigma \mapsto 
\frac{\sigma-1}{\gamma-1}y-(\sigma-1)b \right] \]
defines a map $h^1_{K,V} : \dnrig{}(V)^{\psi=1}_{\Gamma_K} \ra
H^1(K,V)$ which does not depend on the choice of a generator $\gamma$
of $\Gamma_K$, and if $y \in \dfont(V)^{\psi=1} \subset \dnrig{}(V)^{\psi=1}$, 
then $h^1_{K,V}(y)$ coincides with
the cocycle constructed in \cite[I.5]{CC99}.
\end{prop}

\begin{proof}
Our construction closely follows \cite[I.5]{CC99}; to simplify the
proof, we can assume that $\log_p^0(\chi(\gamma))=1$. The fact that if
we start from a different $\gamma$, then the two $h^1_{K,V}$ we get
are the same is left as an easy exercise for the reader. 

We will check that the formula makes sense and that everything is
well-defined.
If $y \in \dnrig{}(V)^{\psi=1}$, then $(\phi-1)y \in
\dnrig{}(V)^{\psi=0}$. By lemma \ref{inv_gam}, there exists $x \in 
\dnrig{}(V)^{\psi=0}$ such that $(\gamma-1)x=(\phi-1)y$. By lemma 
\ref{inv_phi}, there exists $b \in \btrig{}{} \otimes V$ such that
$(\phi-1)b=x$. We then define $h^1_{K,V}(y) \in H^1(K,V)$ 
by the formula:
\[ h^1_{K,V}(y)(\sigma) = 
\frac{\sigma-1}{\gamma-1}y-(\sigma-1)b.\]
Notice that, a priori, $h^1_{K,V}(y) \in H^1(K,\btrig{}{}
\otimes_{\Qp} V)$, but
\begin{align*}
(\phi-1)h^1_{F_n,V}(y)(\sigma) 
& =\frac{\sigma-1}{\gamma-1} (\phi-1)y - (\sigma-1)(\phi-1)b \\
& =\frac{\sigma-1}{\gamma-1} (\gamma-1)x - (\sigma-1)x \\
& =0,
\end{align*}
so that $h^1_{K,V}(y)(\sigma) \in  (\bnrig{}{})^{\phi=1}
\otimes_{\Qp} V = V$. 
In addition, two different choices 
of $b$ differ by an element of 
$(\btrig{}{})^{\phi=1} \otimes_{\Qp} V = V$, 
and therefore give rise
to two cohomologous cocycles.

It is clear that if $y \in \dfont(V)^{\psi=1}
\subset \dnrig{}(V)^{\psi=1}$, then $h^1_{K,V}(y)$ coincides with
the cocycle constructed in \cite[I.5]{CC99}, and that if $y \in
(\gamma-1) \dnrig{}(V)$, then $h^1_{K,V}(y)=0$.
\end{proof}

\section{Explicit formulas for exponential maps}
The goal of this chapter is to give explicit formulas for
Bloch-Kato's maps for a $p$-adic representation $V$, in terms of 
the $(\phi,\Gamma)$-module $\dfont(V)$ attached to $V$. Throughout,
$V$ will be assumed to be crystalline. 

Recall that (cf \cite[III.2]{CC99} or \cite[2.4]{LB02} for example) 
we have maps $\phi^{-n}: \bhol{,F} \ra F_n[[t]]$ which are
characterized by the fact that $\pi$ maps to $\eps^{(n)}e^{t/p^n}-1$.
If $z \in F_n((t)) \otimes_F \dcris(V)$, then the constant 
coefficient (i.e. the
coefficient of $t^0$) of $z$ will be 
denoted by $\partial_V(z) \in \dcris(V)$. 
This notation should not be confused with that for 
the derivation map $\partial$. 

We will make frequent use of
the following fact:
\begin{lemm}\label{tr_psi}
If $y \in \bhol{,F}[1/t] \otimes_F \dcris(V)$, then for any $m \geq n \geq
0$, the element $p^{-m} \on{Tr}_{F_m / F_n} \partial_V(\phi^{-m}(y))
\in F_n \otimes_F \dcris(V)$
does not depend on $m$ and we have:
\[ p^{-m} \on{Tr}_{F_m / F_n} \partial_V(\phi^{-m}(y)) =
\begin{cases}
p^{-n} \partial_V(\phi^{-n}(y)) & \text{if $n \geq 1$} \\
(1-p^{-1}\phi^{-1})\partial_V(y) &  \text{if $n \geq 0$.}
\end{cases} \] 
\end{lemm}

\begin{proof}
Recall that $\phi^{-m}(y)=(1 \otimes
\phi^{-m})y(\eps^{(m)}e^{t/p^m}-1)$, and that $\psi(y)=y$ means that
\[ (1 \otimes \phi)y((1+T)^p-1) = \frac{1}{p} \sum_{\eta^p=1}
y(\eta(1+T)-1). \] 
The lemma then follows from the fact that if $m \geq 2$, then the
conjugates of $\eps^{(m)}$ under $\on{Gal}(F_m/F_{m-1})$ are the $\eta 
\eps^{(m)}$, where $\eta^p=1$, while if $m=1$, then the
conjugates of $\eps^{(1)}$ under $\on{Gal}(F_1/F)$ are the $\eta$, 
where $\eta^p=1$ but $\eta \neq 1$.
\end{proof}

We will henceforth assume that $\log_p(\chi(\gamma_n))=p^n$, so that
$\log_p^0(\chi(\gamma_n))=1$, and in addition
$\nabla_0/(\gamma_n-1)$ acts on $F_n$ by $p^{-n}$.

\Subsection{Bloch-Kato's exponential map}
The goal of this paragraph is to show how to compute Bloch-Kato's map
in terms of the $(\phi,\Gamma)$-module of $V$. 
Let $h \geq 1$ be an integer such that $\on{Fil}^{-h} \dcris(V) =
\dcris(V)$. 
For $i \in \ZZ$, let $\nabla_i$ be the operator 
acting on $F_n((t))$ and $\bhol{,F}$ defined in \ref{notations}.
If $y \in \bhol{,F} \otimes_F \dcris(V)$, then 
the fact that $\on{Fil}^{-h} \dcris(V) =
\dcris(V)$ implies (cf \cite{LBcr} for example) that
$t^h y \in  \dhol{}(V)$, so that if
$y = \sum_{i=0}^d y_i \otimes d_i 
\in (\bhol{,F} \otimes_F \dcris(V))^{\psi=1}$, 
then
\[ \nabla_{h-1} \circ \cdots \circ \nabla_0 
(y) = \sum_{i=0}^d t^h \partial^h y_i \otimes d_i 
\in \dhol(V)^{\psi=1}. \]
One can then apply the operator $h^1_{F_n,V}$ to $\nabla_{h-1} \circ 
\cdots \circ \nabla_0 (y)$, and the main result of this paragraph is:
\begin{theo}\label{exp_bk}
If $y \in (\bhol{,F} \otimes_F \dcris(V))^{\psi=1}$, then
\[  h^1_{F_n,V}
(\nabla_{h-1} \circ  \cdots \circ \nabla_0 (y)) = (-1)^{h-1} (h-1)!
\begin{cases}
\exp_{F_n,V}(p^{-n} \partial_V(\phi^{-n}(y))) & \text{if $n \geq 1$} \\
\exp_{F,V}((1-p^{-1}\phi^{-1})\partial_V(y)) & \text{if $n=0$.}
\end{cases} \]
\end{theo}

\begin{proof}
Because the diagram
\[ \begin{CD} 
F_{n+1} \otimes_F \dcris(V) @>{\exp_{F_{n+1},V}}>> H^1(F_{n+1},V) \\
@V{\on{Tr}_{F_{n+1}/F_n}}VV @V{\on{cor}_{F_{n+1}/F_n}}VV \\
F_n \otimes_F \dcris(V) @>{\exp_{F_n,V}}>> H^1(F_n,V)
\end{CD} \]
is commutative, it is enough to prove the theorem under the further
assumption that $\Gamma_n$ is torsion free. Let us then set 
$y_h=\nabla_{h-1} \circ  \cdots \circ \nabla_0 (y)$.
The cocycle $h^1_{F_n,V}(y_h)$ is defined by 
\[ h^1_{F_n,V}(y_h)(\sigma) = 
\frac{\sigma-1}{\gamma_n-1} y_h - (\sigma-1)b_{n,h} \]
where $b_{n,h}$ is a solution of the equation
$(\gamma_n-1)(\phi-1)b_{n,h} = (\phi-1) y_h$. 
In lemma \ref{inclu_nabla} below, we will prove that:
\[ \nabla_{i-1} \circ \cdots \circ \nabla_1 \circ 
\frac{\nabla_0}{\gamma_n-1} (\bhol{,F})^{\psi=0} 
\subset \left(\frac{t}{\phi^n(\pi)}\right)^i
(\bhol{,F})^{\psi=0}. \]
It is then clear that if one sets
\[ z_{n,h} = \nabla_{h-1} \circ  \cdots \circ \frac{\nabla_0}{\gamma_n-1}
(\phi-1) y, \] then $z_{n,h} \in 
(t/\phi^n(\pi))^h (\bhol{,F})^{\psi=0} \otimes_F \dcris(V) \subset 
\phi^n(\pi^{-h}) 
\dhol(V)^{\psi=0} \subset \dnrig{}(V)^{\psi=0}$.

Let $q_n=\phi^{n-1}(q)=\phi^n(\pi)/\phi^{n-1}(\pi)$;
by lemma \ref{invphiplus} 
(which will be stated and proved below), there exists an element
$b_{n,h} \in \phi^{n-1}(\pi^{-h}) \btrigplus{} \otimes_{\Qp} V$
such that $(\phi-q_n^h)(\phi^{n-1}(\pi^h)b_{n,h})= \phi^n(\pi^h)
z_{n,h}$, so that $(1-\phi)b_{n,h}=z_{n,h}$ with $b_{n,h}
\in \phi^{n-1}(\pi^{-h}) \btrigplus{} \otimes_{\Qp} V$. 

If we set \[ w_{n,h} =  
\nabla_{h-1} \circ  \cdots \circ \frac{\nabla_0}{\gamma_n-1} y, \]
then $w_{n,h}$ and $b_{n,h} \in \bcris \otimes_{\Qp} V$ and
the cocycle $h^1_{F_n,V}(y_h)$ is then given by the formula
$h^1_{F_n,V}(y_h)(\sigma) = (\sigma-1)(w_{n,h} - b_{n,h})$.
Now, $(\phi-1)b_{n,h} = z_{n,h}$ and $(\phi-1)w_{n,h}=z_{n,h}$
as well, so that $w_{n,h} - b_{n,h} \in \bcris^{\phi=1} 
\otimes_{\Qp} V$. 

We can also write $h^1_{F_n,V}(y_h)(\sigma)=(\sigma-1)
(\phi^{-n}(w_{n,h})-\phi^{-n}(b_{n,h}))$. Since we know that $b_{n,h} \in 
\phi^{n-1}(\pi^{-h}) \bcris^+ \otimes_{\Qp} V$, we have
$\phi^{-n}(b_{n,h}) \in \bdR^+ \otimes_{\Qp} V$. By definition of the
Bloch-Kato exponential, the theorem will follow from the fact that:
\[ \phi^{-n}(w_{n,h})- (-1)^{h-1} (h-1)! p^{-n} \partial_V(\phi^{-n}(y)) \in
\bdR^+ \otimes_{\Qp} V. \] 

In order to show this, first notice that
$\phi^{-n}(y)- \partial_V(\phi^{-n}(y)) \in t F_n[[t]] \otimes_F
\dcris(V)$. We can therefore write
$\frac{\nabla_0}{\gamma_n-1} \phi^{-n}(y)
= p^{-n}  \partial_V(\phi^{-n}(y)) + t z_1$ 
and a simple recurrence shows that 
\[ \nabla_{i-1} \circ \cdots \circ
\frac{\nabla_0}{1-\gamma_n} \phi^{-n}(y)
= (-1)^{i-1} (i-1)! p^{-n} \partial_V(\phi^{-n}(y)) + t^i z_i, \] 
with $z_i \in F_n[[t]] \otimes_F
\dcris(V)$. By taking $i=h$, we see that
$\phi^{-n}(w_{n,h}) - (-1)^{h-1} (h-1)! p^{-n} \partial_V(\phi^{-n}(y)) \in
\bdR^+ \otimes_{\Qp} V$, since we chose $h$ such that
$t^h \dcris(V) \subset \bdR^+
\otimes_{\Qp} V$.
\end{proof}
 
We will now prove the following two technical lemmas which were used above:

\begin{lemm}\label{inclu_nabla}
If $n \geq 1$, then
$\nabla_0/(\gamma_n-1) (\bhol{,F})^{\psi=0} \subset (t/\phi^n(\pi))
(\bhol{,F})^{\psi=0}$ so that if $i \geq 1$, then: 
\[ \nabla_{i-1} \circ \cdots \circ \nabla_1 \circ 
\frac{\nabla_0}{\gamma_n-1} (\bhol{,F})^{\psi=0}
\subset \left(\frac{t}{\phi^n(\pi)}\right)^i
(\bhol{,F})^{\psi=0}. \]
\end{lemm}

\begin{proof}
Since $\nabla_i = t \cdot d/dt -i$, the second claim follows easily
from the first one. By the standard properties of $p$-adic holomorphic
functions, what we need to do is to show that if
$x \in (\bhol{,F})^{\psi=0}$, then
$(\nabla_0/(\gamma_n-1) x)(\eps^{(m)}-1)=0$ for all
$m \geq n+1$. 

On the one hand,
up to a scalar factor, one has for $m \geq n+1$:
$(\nabla_0/(\gamma_n-1) x)(\eps^{(m)}-1) = \on{Tr}_{F_m/F_n}
x(\eps^{(m)}-1)$ as can be seen from the fact that  
$\nabla_0/(\gamma_n-1) = \lim_{\Gamma \ni \eta \ra 1}
(\eta-1)/(\gamma_n-1) \cdot \log_p^{-1}(\chi(\eta))$.
On the other hand, the fact that $\psi(x)=0$ implies that for every 
$m \geq 2$, $\on{Tr}_{F_m/F_{m-1}} x(\eps^{(m)}-1)=0$. This completes
the proof.
\end{proof}

\begin{lemm}\label{invphiplus}
If $\alpha \in \btrigplus{}$, then there exists $\beta \in
\btrigplus{}$ such that $(\phi-q_n^h)\beta=\alpha$.  
\end{lemm}

\begin{proof}
By \cite[prop 2.19]{LB02} applied to the case $r=0$, the ring
$\btplus$ is dense in $\btrigplus{}$ for the Fr\'echet
topology. Hence, if $\alpha \in \btrigplus{}$, then there exists
$\alpha_0 \in \btplus$ such that $\alpha - \alpha_0 = \phi^n(\pi^h) \alpha_1$
with $\alpha_1 \in \btrigplus{}$. The map $\phi-q_n^h : \btplus \ra
\btplus$ is surjective, because  $\phi-q_n^h : \atplus \ra
\atplus$ is surjective, as can be seen by reducing modulo $p$. One can
therefore write $\alpha_0 = (\phi-q_n^h)\beta_0$. Finally (see the
proof of lemma \ref{inv_phi}), there exists $\beta_1 \in \btrigplus{}$
such that $\alpha_1=(\phi-1)\beta_1$, so that 
$\phi^n(\pi^h) \alpha_1 = (\phi-q_n^h) (\phi^{n-1}(\pi^h) \beta_1)$.
\end{proof}

\Subsection{Bloch-Kato's dual exponential map}
In the previous paragraph, we showed how to compute Bloch-Kato's
exponential map for $V$. We shall now do the same for the dual
exponential map. The starting point is Kato's formula, which we recall
below:
\begin{prop}\label{kato_formula}
If $V$ is a de Rham representation, then the map from $\ddR(V)$ to
$H^1(F,\bdR \otimes_{\Qp} V)$ defined by 
$x \mapsto \left[ g \mapsto \log(\chi(\overline{g}))x \right]$  
is an isomorphism, and the dual
exponential map $\exp^*_{V^*(1)} : H^1(F,V) \ra \ddR(V)$
is equal to the composition of the map $H^1(F,V) \ra
H^1(F,\bdR \otimes_{\Qp} V)$ with the inverse of this isomorphism.
\end{prop}

Let us point out that the image of $\exp^*_{V^*(1)}$ is included in
$\on{Fil}^0 \ddR(V)$ and that its kernel is $H^1_g(F,V)$. Let us
return to a crystalline representation $V$ of $G_F$. 
We then have the following formula,
which is to be found in \cite[IV.2.1]{CC99}:

\begin{theo}\label{exp_duale}
If $y \in \dnrig{}(V)^{\psi=1}$, then
\[ \exp^*_{F_n, V^*(1)}(h^1_{F_n,V}(y)) = 
\begin{cases}
p^{-n} \partial_V(\phi^{-n}(y)) & \text{if $n \geq 1$} \\
(1-p^{-1}\phi^{-1})\partial_V(y) &  \text{if $n = 0$.}
\end{cases} \]
\end{theo}

\begin{proof}
Since the following diagram 
\[ \begin{CD} 
H^1(F_{n+1},V) @>{\exp^*_{F_{n+1},V^*(1)}}>> F_{n+1} \otimes_F
\dcris(V) \\ 
@V{\on{cor}_{F_{n+1}/F_n}}VV @V{\on{Tr}_{F_{n+1}/F_n}}VV \\
H^1(F_n,V) @>{\exp^*_{F_n,V^*(1)}}>> F_n \otimes_F \dcris(V)
\end{CD} \]
is commutative, we only need to prove the theorem when $\Gamma_n$ is
torsion free. We then have \[ h^1_{F_n,V}(y)(\sigma) 
= \frac{\sigma-1}{\gamma_n-1} y - (\sigma - 1)
b, \] where $(\gamma_n-1)(\phi-1)b=(\phi-1)y$. 
Since $b \in \btrig{}{}$, there exists
$m \gg 0$ such that $b \in \btrig{,r_m}{} \otimes_{\Qp} V$.
Recall that (cf \cite[2.4]{LB02}) the map $\phi^{-m}$ embeds
$\btrig{,r_m}{}$ into $\bdR^+$.
We can then write
\[ h^1(y)(\sigma) = \frac{\sigma-1}{\gamma_n-1} \phi^{-m}(y) - (\sigma - 1)
\phi^{-m}(b), \] and $\phi^{-m}(b) \in \bdR^+ \otimes_{\Qp} V$. 
In addition,
$\phi^{-m}(y) \in F_m((t)) \otimes_F \dcris(V)$ and $\gamma_n-1$
is invertible on $t^k F_m \otimes_F \dcris(V)$ for every $k \neq
0$. This shows that the cocycle
$h^1_{F_n,V}(y)$ is cohomologous in $H^1(F_n,\bdR
\otimes_{\Qp} V)$ to
\[ \sigma \mapsto \frac{\sigma-1}{\gamma_n-1}(\partial_V \phi^{-m}(y))
\] which is itself cohomologous 
(since $\gamma_n-1$ is invertible on $F_m^{\on{Tr}_{F_m/F_n}=0}$)
to
\begin{equation*}
\sigma \mapsto \frac{\sigma-1}{\gamma_n-1}\left(p^{n-m}
\on{Tr}_{F_m/F_n}
\partial_V (\phi^{-m}(y))\right) 
 = \sigma \mapsto p^{-n} \log(\chi(\overline{\sigma}))
p^{n-m} \on{Tr}_{F_m/F_n}
\partial_V (\phi^{-m}(y)). 
\end{equation*}
It follows from this and Kato's formula (proposition
\ref{kato_formula}) that
\begin{equation*} \exp^*_{F_n, V^*(1)}(h^1_{F_n,V}(y))=p^{-m} \on{Tr}_{F_m/F_n}
\partial_V (\phi^{-m}(y)) 
= \begin{cases}
p^{-n} \partial_V(\phi^{-n}(y)) & \text{if $n \geq 1$} \\
(1-p^{-1}\phi^{-1})\partial_V(y) &  \text{if $n = 0$.}
\end{cases} \end{equation*}
\end{proof}

\Subsection{Iwasawa theory of $p$-adic representations}
In this specific paragraph, $V$ 
can be taken to be an arbitrary representation of 
$G_K$. Recall that the Iwasawa cohomology groups $H^i_{Iw}(K,V)$ are
defined by $H^i_{Iw}(K,V) = \Qp \otimes_{\Zp} H^i_{Iw}(K,T)$ where $T$
is any $G_K$-stable lattice of $V$, and \[ H^i_{Iw}(K,T) = 
\projlim_{\on{cor}_{K_{n+1}/K_n}} 
H^i(K_n,T). \]
Each of the $H^i(K_n,T)$ is a $\Zp[\Gamma/\Gamma_n]$-module, and 
$H^i_{Iw}(K_n,T)$ is then endowed with the structure of a
$\Lambda$-module where $\Lambda=\Zp[[\Gamma]]=\Zp[\Delta] \otimes_{\Zp}
\Zp[[\Gamma_1]]$.  The $H^i_{Iw}(K,V)$ have been studied in detail by
Perrin-Riou, who proved the following:
\begin{prop}
If $V$ is a $p$-adic representation of $G_K$, then $H^i_{Iw}(K,V)=0$
whenever $i \neq 1$, $2$. In addition:
\begin{enumerate}
\item the torsion sub-module of $H^1_{Iw}(K,V)$ is a 
$\Qp \otimes_{\Zp} \Lambda$-module
isomorphic to $V^{H_K}$, and then $H^1_{Iw}(K,V)/V^{H_K}$ is a free
$\Qp \otimes_{\Zp} \Lambda$-module of rank $[K:\Qp]d$;
\item $H^2_{Iw}(K,V) = V(-1)^{H_K}$.
\end{enumerate}
\end{prop}

If $y \in \dfont(T)^{\psi=1}$, then the sequence of $\{h^1_{F_n,V}(y)\}_n$
is compatible for the corestriction maps, and therefore defines an
element of $H^1_{Iw}(K,T)$. The following theorem is due to
Fontaine and is proved in \cite[II.1]{CC99}:
\begin{theo}\label{dvpsi}
The map $y \mapsto \projlim_n h^1_{K_n,V}(y)$ defines an isomorphism from 
$\dfont(T)^{\psi=1}$ to $H^1_{Iw}(K,T)$ and from 
$\dfont(V)^{\psi=1}$ to $H^1_{Iw}(K,V)$.
\end{theo}

Notice that $V^{H_K} \subset \dfont(V)^{\psi=1}$, and it 
its $\Qp \otimes_{\Zp} \Lambda$-torsion submodule.
In addition, it is shown in \cite[II.3]{CC99} that
the modules $\dfont(V)/(\psi-1)$ and $H^2_{Iw}(K,V)$ are naturally
isomorphic. One can nicely summarize the results 
of this paragraph as follows:

\begin{coro}\label{compiw}
The complex of $\Qp \otimes_{\Zp} \Lambda$-modules \[ \begin{CD} 0 
@>>> \dfont(V) @>{1-\psi}>> \dfont(V) @>>> 0 \end{CD} \]
computes the Iwasawa cohomology of $V$. 
\end{coro}

\Subsection{Perrin-Riou's exponential map}
By using the results of the previous paragraphs, we can give a
``uniform'' formula  for the image of $y \in (\bhol{,F} 
\otimes_F \dcris(V))^{\psi=1}$ in $H^1(F_n,V(j))$ under the
composition of the following maps:
\begin{multline*} 
\begin{CD} \left( \bhol{,F} \otimes_F \dcris(V) \right)^{\psi=1}
@>{\nabla_{h-1} \circ \cdots \circ \nabla_0}>> 
\dnrig{}(V)^{\psi=1} @>{\otimes e_j}>> \end{CD} \\
\begin{CD}
\dnrig{}(V(j))^{\psi=1} @>{h^1_{F_n,V(j)}}>>
H^1(F_n,V(j)). \end{CD} 
\end{multline*}

\begin{theo}\label{recip}
If $y \in (\bhol{,F} \otimes_F \dcris(V))^{\psi=1}$, 
and $h \geq 1$ is such that
$\on{Fil}^{-h}\dcris(V)=\dcris(V)$,
then for all $j$ with $h+j \geq 1$, we have:
\begin{multline*}
h^1_{F_n,V(j)}(\nabla_{h-1} \circ \cdots \circ \nabla_0 (y) \otimes
e_j)=
\\ 
(-1)^{h+j-1} (h+j-1)!
\begin{cases}
\exp_{F_n,V(j)}(p^{-n} \partial_{V(j)}(\phi^{-n} 
(\partial^{-j}y \otimes t^{-j}e_j)))
& \text{if $n \geq 1$} \\
\exp_{F,V(j)}((1-p^{-1} \phi^{-1})\partial_{V(j)} 
(\partial^{-j}y \otimes t^{-j}e_j))
& \text{if $n=0$,}
\end{cases} 
\end{multline*}
while if $h+j \leq 0$, then we have:
\begin{multline*}
\exp^*_{F_n,V^*(1-j)}(h^1_{F_n,V(j)}(\nabla_{h-1} \circ 
\cdots \circ \nabla_0 (y) \otimes
e_j))=
\\
(-h-j)!^{-1}
\begin{cases} 
p^{-n} \partial_{V(j)}(\phi^{-n} 
(\partial^{-j}y \otimes t^{-j}e_j))
& \text{if $n \geq 1$} \\
(1-p^{-1} \phi^{-1})\partial_{V(j)} 
(\partial^{-j}y \otimes t^{-j}e_j)
& \text{if $n = 0$.}
\end{cases}
\end{multline*}
\end{theo}

\begin{proof}
If $h+j \geq 1$, then the following diagram is commutative:
\[ \begin{CD}
\dhol(V)^{\psi=1} @>{\otimes e_j}>> \dhol(V(j))^{\psi=1}  \\
@A{\nabla_{h-1} \circ \cdots \circ \nabla_0}AA 
@A{\nabla_{h+j-1} \circ \cdots \circ \nabla_0}AA \\
\left( \bhol{,F} \otimes_F \dcris(V) \right)^{\psi=1}
@>{\partial^{-j} \otimes t^{-j} e_j}>>
\left( \bhol{,F} \otimes_F \dcris(V(j)) \right)^{\psi=1}.
\end{CD} \]
and the theorem is then a straightforward consequence of theorem 
\ref{exp_bk} applied to $\partial^{-j}y \otimes t^{-j}e_j$, $h+j$
and $V(j)$ (which are the $j$-th twists of $y$, $h$ and $V$).

If on the other hand $h+j \leq 0$, and $\Gamma_n$ is torsion free, 
then theorem \ref{exp_duale} shows that
\begin{equation*} 
\exp^*_{F_n,V^*(1-j)}
(h^1_{F_n,V(j)}(\nabla_{h-1} \circ \cdots 
\circ \nabla_0 (y) \otimes e_j))  = 
p^{-n} \partial_{V(j)}(\phi^{-n}(\nabla_{h-1} \circ \cdots 
\circ \nabla_0 (y) \otimes e_j)) 
\end{equation*} 
in $\dcris(V(j))$, 
and a short computation involving Taylor series shows that 
\begin{equation*}
p^{-n} \partial_{V(j)}(\phi^{-n}(\nabla_{h-1} \circ \cdots 
\circ \nabla_0 (y) \otimes e_j)) = 
(-h-j)!^{-1} p^{-n} \partial_{V(j)}(\phi^{-n} 
(\partial^{-j}y \otimes t^{-j}e_j)). 
\end{equation*}
Finally, to get the case $n=0$, one just needs to corestrict.
\end{proof}

We shall now use the above result to construct Perrin-Riou's
exponential map. One has an exact sequence
\begin{multline*}
\begin{CD} 0 @>>> \oplus_{k=0}^h t^k \dcris(V)^{\phi=p^{-k}} @>>> 
\left( \bhol{,F} \otimes_F \dcris(V) \right)^{\psi=1} 
@>{1-\phi}>> \end{CD}
\\ \begin{CD}
(\bhol{,F})^{\psi=0} \otimes_F \dcris(V)
@>{\Delta}>> \oplus_{k=0}^h \left(\frac{\dcris(V)}{1-p^k 
\phi}\right)(k) @>>> 0, \end{CD}
\end{multline*}
where $\Delta(f)=\oplus_{k=0}^h \partial^k(f)(0)$. 
If $f \in ((\bhol{,F})^{\psi=0} \otimes_F \dcris(V))^{\Delta=0}$, 
then there exists $y \in ( \bhol{,F} \otimes_F \dcris(V) )^{\psi=1}$ 
such that $f = (1-\phi)y$, and 
$\nabla_{h-1} \circ \cdots \circ \nabla_0 (y)$ does not depend upon
the choice of such a $y$ (unless $\Qp(h) \subset V$), and one deduces
from this a well-defined map (if $\Qp(k) \subset V$, we therefore
require that $h \geq k+1$):
\[ \Omega_{V,h}:
((\bhol{,F})^{\psi=0} \otimes_F \dcris(V))^{\Delta=0} \ra
\dhol(V)^{\psi=1}, \]
given by
$\Omega_{V,h}(f) = \nabla_{h-1} \circ \cdots \circ \nabla_0 (y)$.

\begin{theo}\label{bpr_const}
If $V$ is a crystalline representation and $h \geq 1 $ is such that we have
$\on{Fil}^{-h} \dcris(V) = \dcris(V)$, then the map 
\[ \Omega_{V,h} : ((\bhol{,F})^{\psi=0} \otimes_F \dcris(V))^{\Delta=0} 
\ra \dhol(V)^{\psi=1} / V^{H_F} \] 
which takes $x \in ((\bhol{,F})^{\psi=0} \otimes_F
\dcris(V))^{\Delta=0}$ to $\nabla_{h-1} \circ \cdots \circ \nabla_0
((1-\phi)^{-1}(x))$ 
is well-defined and coincides
with Perrin-Riou's exponential map.
\end{theo}

\begin{proof}
The map $\Omega_{V,h}$ is well defined because the kernel of $1-\phi$
is killed by $\nabla_{h-1} \circ \cdots \circ \nabla_0$, except for
$t^h \otimes (t^{-h}e_h)$, which is mapped to $\Qp(h) \subset V^{H_F}$.

The fact that $\Omega_{V,h}$  coincides
with Perrin-Riou's exponential map
follows directly from theorem \ref{recip} above applied to those 
$j$'s for which $h+j \geq 1$, compared with \cite[3.2.3]{BP94} (see
remark \ref{crazy} however).
\end{proof}

\begin{rema}\label{no_kill}
By the above remarks, if 
$V$ is a crystalline representation and $h \geq 1 $ is such that we have
$\on{Fil}^{-h} \dcris(V) = \dcris(V)$ and $\Qp(h) \not \subset V$, 
then the map 
\[ \Omega_{V,h} : ((\bhol{,F})^{\psi=0} \otimes_F \dcris(V))^{\Delta=0} 
\ra \dhol(V)^{\psi=1} \] 
which takes $x \in ((\bhol{,F})^{\psi=0} \otimes_F
\dcris(V))^{\Delta=0}$ to $\nabla_{h-1} \circ \cdots \circ \nabla_0
((1-\phi)^{-1}(x))$ 
is well-defined, without having to kill the $\Lambda$-torsion 
of $H^1_{Iw}(F,V)$.
\end{rema}

\begin{rema}\label{twist_bpr}
It is clear from theorem \ref{recip} that we have:
\[ \Omega_{V,h}(x) \otimes e_j = \Omega_{V(j),h+j}(\partial^{-j}x
\otimes t^{-j}e_j)
\quad\text{and}\quad
\Omega_{V,h}(\nabla_h (x)) = \Omega_{V,h+1}(x), \]
and following Perrin-Riou, one can use these formulas to extend
the definition of $\Omega_{V,h}$ to all $h \in \ZZ$.
\end{rema}

\Subsection{The explicit reciprocity formula}
In this paragraph, we shall recall Perrin-Riou's explicit
reciprocity formula. There is a map $\mathcal{H}(\Gamma_F) 
\ra (\bhol{,\Qp})^{\psi=0}$
which sends $f(\gamma)$ to $f(\gamma)(1+\pi)$. This map is a bijection
and its inverse is the Mellin transform so that 
if $g(\pi) \in (\bhol{,\Qp})^{\psi=0}$, then
$g(\pi)=\on{Mel}(g)(1+\pi)$. See \cite[B.2.8]{BP00} for a reference, where 
Perrin-Riou has also extended $\on{Mel}$ to $(\bnrig{}{,\Qp})^{\psi=0}$.
If $f,g \in \bnrig{}{,\Qp}$ then we define $f*g$ by the formula 
$\on{Mel}(f*g)=\on{Mel}(f)\on{Mel}(g)$. 
Let $[-1] \in \Gamma_F$ be the element such that $\chi([-1])=-1$, and
let $\iota$ be the involution of $\Gamma_F$ which sends $\gamma$ to
$\gamma^{-1}$. The operator $\partial^j$ on $(\bhol{,\Qp})^{\psi=0}$
corresponds to $\on{Tw}_j$ on $\Gamma_F$. For instance, it is a bijection.
We will make use of the facts that $\iota \circ
\partial^j = \partial^{-j} \circ \iota$ and that 
$[-1] \circ \partial^j = (-1)^j \partial^j \circ [-1]$.

If $V$ is a crystalline representation, then
the natural maps \[ \begin{CD}
\dcris(V) \otimes_F \dcris(V^*(1)) 
@>>> \dcris(V \otimes_{\Qp} V^*(1)) @>>> \dcris(\Qp(1))
@>{\on{Tr_{F/\Qp}}}>> \Qp \end{CD} \] allow us to define a perfect
pairing $[\cdot,\cdot]_V: \dcris(V) \times \dcris(V^*(1))$ which we
extend by linearity to 
\[ [\cdot,\cdot]_V : (\bhol{,F} \otimes_F \dcris(V))^{\psi=0} \times 
(\bhol{,F} \otimes_F \dcris(V^*(1)))^{\psi=0} \ra
(\bhol{,\Qp})^{\psi=0} \]
by the formula $[f(\pi) \otimes d_1, g(\pi) \otimes
d_2]_V = (f*g)(\pi) [d_1,d_2]_V$.

We can also define a semi-linear (with respect to $\iota$)
pairing \[ \langle \cdot, \cdot \rangle_V : 
\dhol(V)^{\psi=1} \times \dhol(V^*(1))^{\psi=1} \ra
(\bhol{,\Qp})^{\psi=0} \] by the formula
\[ \langle y_1, y_2 \rangle_V
= \projlim_n \sum_{\tau \in \Gamma_F/\Gamma_n} \langle
\tau^{-1}(h^1_{F_n,V}(y_1)),h^1_{F_n,V^*(1)}(y_2) 
\rangle_{F_n,V} \cdot \tau(1+\pi) \]
where the pairing $\langle \cdot, \cdot \rangle_{F_n,V}$ is given by the cup
product: \[ \langle \cdot, \cdot \rangle_{F_n,V} : H^1(F_n,V) \times
H^1(F_n,V^*(1)) \ra H^2(F_n, \Qp(1)) \simeq \Qp. \]
The pairing $\langle \cdot, \cdot \rangle_V$ 
satisfies the relation $\langle \gamma_1 x_1, \gamma_2 x_2 \rangle_V
=\gamma_1 \iota(\gamma_2) \langle x_1, x_2 \rangle_V$.
Perrin-Riou's explicit reciprocity formula 
(proved by Colmez \cite{Co98}, Benois 
\cite{Be00} and Kato-Kurihara-Tsuji (unpublished))
is then:
\begin{theo}\label{rec_formula}
If $x_1 \in (\bhol{,F} \otimes_F \dcris(V))^{\psi=0}$ and 
$x_2 \in (\bhol{,F} \otimes_F \dcris(V^*(1)))^{\psi=0}$, then for
every $h$, we have:
\[ (-1)^h  \langle \Omega_{V,h}(x_1), [-1] \cdot
\Omega_{V^*(1),1-h}(x_2) \rangle_V = - [x_1,\iota(x_2)]_V. \]
\end{theo}

\begin{proof}
By the theory of $p$-adic interpolation, it is enough to prove that
if $x_i=(1-\phi)y_i$ with 
$y_1 \in (\bhol{,F} \otimes_F \dcris(V))^{\psi=1}$ and 
$y_2 \in (\bhol{,F} \otimes_F \dcris(V^*(1)))^{\psi=1}$ then for all
$j \gg 0$:
\[ \left(\partial^{-j}  (-1)^h  \langle \Omega_{V,h}(x_1), [-1] \cdot
\Omega_{V^*(1),1-h}(x_2) \rangle_V \right)(0) = - \left(\partial^{-j}
[x_1,\iota(x_2)]_V\right)(0). \] 

The above formula is equivalent to (a):
\begin{multline*}
  (-1)^h  \langle h^1_{F,V(j)} \Omega_{V(j),h+j} 
(\partial^{-j} x_1 \otimes t^{-j}e_j), (-1)^j 
 h^1_{F,V^*(1-j)} \Omega_{V^*(1-j),1-h-j}
(\partial^j x_2 \otimes t^j e_{-j}) \rangle_{F,V(j)} \\  = 
- [\partial_{V(j)} (\partial^{-j} x_1 \otimes t^{-j}e_j),
\partial_{V^*(1-j)}(\partial^j x_2 \otimes t^j e_{-j})]_{V(j)}. 
\end{multline*}
By combining theorems \ref{recip} and \ref{bpr_const}
with remark \ref{twist_bpr} 
we see that for $j \gg 0$: 
\begin{multline*} h^1_{F,V(j)} \Omega_{V(j),h+j} 
(\partial^{-j} x_1 \otimes t^{-j}e_j) \\ = 
(-1)^{h+j-1} \exp_{F,V(j)}((h+j-1)!(1-p^{-1} \phi^{-1})\partial_{V(j)} 
(\partial^{-j}y_1 \otimes t^{-j}e_j)), \end{multline*}
and that
\begin{multline*}   h^1_{F,V^*(1-j)} \Omega_{V^*(1-j),1-h-j}
(\partial^j x_2 \otimes t^j e_{-j}) \\ = 
(\exp^*_{F,V^*(1-j)})^{-1} (h+j-1)!^{-1}((1-p^{-1} \phi^{-1})
\partial_{V^*(1-j)} 
(\partial^j y_2 \otimes t^j e_{-j})).\end{multline*}
Using the fact that by definition,
if $x \in \dcris(V(j))$ and $y \in H^1(F,V(j))$ then 
\[ [x,\exp^*_{F,V^*(1-j)} y]_{V(j)}
=\langle \exp_{F,V(j)} x , y \rangle_{F,V(j)}, \]
we see that (b):
\begin{multline*}
\langle h^1_{F,V(j)} \Omega_{V(j),h+j} 
(\partial^{-j} x_1 \otimes t^{-j}e_j),
 h^1_{F,V^*(1-j)} \Omega_{V^*(1-j),1-h-j}
(\partial^j x_2 \otimes t^j e_{-j}) \rangle_{F,V(j)} \\  = 
(-1)^{h+j-1} 
[(1-p^{-1}\phi^{-1})\partial_{V(j)} (\partial^{-j} y_1 \otimes t^{-j}e_j),
(1-p^{-1}\phi^{-1})
\partial_{V^*(1-j)}(\partial^j y_2 \otimes t^j e_{-j})]_{V(j)}. 
\end{multline*}
It is easy to see that under $[\cdot,\cdot]$, the adjoint of
$(1-p^{-1}\phi^{-1})$ is $1-\phi$, and that if $x_i=(1-\phi)y_i$, then
\begin{align*} \partial_{V(j)} (\partial^{-j} x_1 \otimes t^{-j}e_j) 
& = (1-\phi) \partial_{V(j)} (\partial^{-j} y_1 \otimes t^{-j}e_j), 
\\
\partial_{V^*(1-j)} (\partial^j x_2 \otimes t^j e_{-j}) 
& = (1-\phi) \partial_{V^*(1-j)} (\partial^j y_2 \otimes t^j e_{-j}), 
\end{align*} so
that (b) implies (a), and this proves the
theorem.
\end{proof}

\begin{rema}\label{crazy}
One should be careful with all the signs involved in those
formulas. Perrin-Riou has changed the definition 
of the $\ell_i$ operators from
\cite{BP94} to \cite{BP99}
(the new $\ell_i$ is minus the old $\ell_i$). 
The reciprocity formula which is stated
in \cite[4.2.3]{BP99} does not seem (to me) 
to have the correct sign. On the other
hand, the formulas of \cite{Be00,Co98} do seem to give the correct signs, 
but one should be careful that \cite[IX.4.5]{Co98} uses a different 
definition for one of the pairings,
and that the signs in \cite[IV.3.1]{CC99} and \cite[VII.1.1]{Co98}
disagree. Our definitions of $\Omega_{V,h}$ and of the pairing agree 
with Perrin-Riou's ones (as they are given in \cite{BP99}).
\end{rema}

\newcommand{\appendixname}{Appendix}

\appendix\section{The structure of $\dfont(T)^{\psi=1}$}
The goal of this paragraph is to prove a theorem which says that 
for a crystalline representation $V$, $\dfont(V)^{\psi=1}$ is quite ``small''.
See theorem \ref{quite_small} for a precise statement.

Let $V$ be a
crystalline representation of $G_F$
and let $T$ denote a
$G_F$-stable lattice of $V$.
The following proposition, which improves slightly upon the results of
N. Wach \cite{Wa96}, is proved in detail in \cite {LBcr}:

\begin{prop}\label{const_nt}
If $T$ is a lattice in a positive crystalline 
representation $V$, then there exists a unique
sub $\aplus_F$-module $N(T)$ of $\dfont^+(T)$,
which satisfies the following conditions:
\begin{enumerate}
\item $N(T)$ is free of rank $d$;
\item $N(T)$ is stable under
the action of $\Gamma_F$; 
\item the action  of $\Gamma_F$ is trivial on $N(T) /\pi N(T)$;
\item there exists an integer $s \geq 0$ such 
that $\pi^s \dfont^+(T) \subset N(T)$. 
\end{enumerate}
In this case, $N(T)$ is stable by $\phi$, 
and the sub $\bplus_F$-module $N(V) = N(T) \otimes_{\aplus_F} \bplus_F$ of
$\dfont^+(V)$ satisfies the corresponding conditions.
\end{prop}

Notice that $N(T(-1))=\pi N(T) \otimes e_{-1}$.
When $V$ is no longer positive, we therefore
define $N(T)$ as $\pi^{-h}
N(T(-h)) \otimes e_h$, for $h$ large enough so that $V(-h)$ is positive.

\begin{prop}
If $T$ is a lattice in a crystalline representation $V$ of $G_F$,
whose Hodge-Tate weights are in $[a;b]$, then
$N(T)$ is the unique sub-$\aplus_F$-module of
$\dfont^+(T)[1/\pi]$ which is free of rank $d$, stable by $\Gamma_F$
with the action of $\Gamma_F$ trivial on $N(T) /\pi N(T)$, and such
that $N(T)[1/\pi]=\dfont^+(T)[1/\pi]$. 

In addition, we have $\phi(\pi^b N(T)) \subset \pi^b N(T)$
and $\pi^b N(T) / \phi^*(\pi^b N(T))$ is killed by $q^{b-a}$.
The functor $T \mapsto N(T)$ is an equivalence of categories between
the category of lattices in crystalline representations and the
category of modules satisfying all the above conditions.
\end{prop}

We shall now show that $\dfont(V)^{\psi=1}$ is not very far from being
included in $N(V)$. Indeed:

\begin{theo}\label{quite_small}
If $V$ is a crystalline representation of $G_F$,
whose Hodge-Tate weights are in $[a;b]$, then
$\dfont(V)^{\psi=1} \subset \pi^{a-1} N(V)$. 

If in addition $V$ has no quotient isomorphic to $\Qp(a)$, then 
$\dfont(V)^{\psi=1} \subset \pi^a N(V)$. 
\end{theo}

Before we prove the above statement, we will need a few results
concerning the action of $\psi$ on $\dfont(T)$. In lemmas
\ref{first_lem} through \ref{last_lem}, we will assume that the
Hodge-Tate weights of $V$ are $\geq 0$. In particular, 
$N(T) \subset \phi^* N(T)$ so that $\psi(N(T)) \subset N(T)$. 

\begin{lemm}\label{exo}
If $m \geq 1$, then there exists a polynomial $Q_m(X) \in \Zp[X]$ such that
$\psi(\pi^{-m}) = \pi^{-m}(p^{m-1}+\pi Q_m(\pi))$.
\end{lemm}

\begin{proof}
By the definition of $\psi$, it is enough to show that if  
$m \geq 1$, there exists a polynomial $Q_m(X) \in \ZZ[X]$ such that
\[ \frac{1}{p} \sum_{\eta^p=1} \frac{1}{(\eta (1+X)-1)^m} = 
\frac{p^{m-1}+((1+X)^p-1)Q_m((1+X)^p-1)}{((1+X)^p-1)^m}, \]
which is left to the reader.
\end{proof}

\begin{lemm}\label{first_lem}
If $k \geq 1$, then $\psi(p \dfont(T) + \pi^{-(k+1)}N(T)) \subset
p \dfont(T) + \pi^{-k}N(T)$.  In addition, 
 $\psi(p \dfont(T) + \pi^{-1}N(T)) \subset
p \dfont(T) + \pi^{-1}N(T)$
\end{lemm}

\begin{proof}
If $x \in N(T)$, then one can write $x=\sum \lambda_i \phi(x_i)$ with
$\lambda_i \in \aplus_F$ and $x_i \in N(T)$, so that 
$\psi(\pi^{-(k+1)} x) = \sum \psi(\pi^{-(k+1)} \lambda_i) x_i$. By the
preceding lemma, $\psi(\pi^{-(k+1)} \lambda_i) \in p \mathbf{A}_F + 
\pi^{-k}\aplus_F$ whenever $k \geq 1$. 
The lemma follows easily, and the second claim is proved in
the same way.
\end{proof}

\begin{lemm}\label{mid_lem}
If $k \geq 1$ and $x \in \dfont(T)$ 
has the property that $\psi(x)-x \in p \dfont(T) + \pi^{-k}N(T)$,
then $x \in p \dfont(T) + \pi^{-k}N(T)$.
\end{lemm}

\begin{proof}
Let $\ell$ be the smallest integer $\geq 0$ such 
that $x \in p \dfont(T) + \pi^{-\ell}N(T)$. If $\ell \leq k$, then we are
done and otherwise lemma \ref{first_lem} shows that
$\psi(x) \in p \dfont(T) + \pi^{-(\ell-1)}N(T)$, so that $\psi(x)-x$ would be
in $p \dfont(T) + \pi^{-\ell}N(T)$ but not $p \dfont(T) +
\pi^{-(\ell-1)}N(T)$, a contradiction if $\ell > k$. 
\end{proof}

\begin{lemm}\label{last_lem}
We have $\dfont(T)^{\psi=1} \subset \pi^{-1}N(T)$.
\end{lemm}

\begin{proof}
We shall prove by induction that $\dfont(T)^{\psi=1} \subset 
p^k \dfont(T) + \pi^{-1}N(T)$ for $k \geq 1$. 
Let us start with the case $k=1$.
If $x \in \dfont(T)^{\psi=1}$,
then there exists some $j \geq 1$ such that $x \in 
p \dfont(T) + \pi^{-j}N(T)$. If $j=1$ we are done and otherwise
the fact that $\psi(x)=x$ combined with lemma \ref{first_lem} shows
that $j$ can be decreased by $1$. This proves our claim for $k=1$.

We will now assume our claim 
to be true for $k$ and prove it for $k+1$. If $x
\in \dfont(T)^{\psi=1}$, we can therefore write $x=p^k y + n$ where
$y \in \dfont(T)$ and $n \in \pi^{-1}N(T)$. Since $\psi(x)=x$, we
have $\psi(n)-n = p^k (\psi(y)-y)$ so that $\psi(y)-y \in \pi^{-1}N(T)$
(this is because $p^k \dfont(T) \cap N(T) = p^k N(T)$). By lemma
\ref{mid_lem}, this implies that $y \in p \dfont(T) + \pi^{-1}N(T)$,
so that we can write $x=p^k(py'+n')+n=p^{k+1}y'+(p^k n'+n)$, 
and this proves our claim.

Finally, it is clear that our claim implies the lemma: if one can
write $x = p^k y_k + n_k$, then the $n_k$ will converge for the
$p$-adic topology to a $n \in \pi^{-1} N$ such that $x=n$.  
\end{proof}

\begin{proof}[Proof of theorem \ref{quite_small}]
Clearly, it is enough to show that if $T$ is a $G_F$-stable lattice of
$V$, then $\dfont(T)^{\psi=1} \subset \pi^{a-1} N(T)$.
It is also
clear that we can twist $V$ as we wish, and we shall now assume that
the Hodge-Tate weights of $V$ are in $[0;h]$. In this case, the
theorem says that $\dfont(T)^{\psi=1} \subset \pi^{-1} N(T)$, which is the
content of lemma \ref{last_lem} above.

Let us now prove that if a positive $V$ has no quotient isomorphic
to $\Qp$, then actually $\dfont(T)^{\psi=1} \subset N(T)$.
Recall that $N(T) 
\subset \phi^*(N(T))$, since the Hodge-Tate weights of $V$ are $\geq 0$,
so that if $e_1, \cdots, e_d$ is a
basis of $N(T)$, then there exists $q_{ij} \in \aplus_F$ such that
$e_i = \sum_{j=1}^d q_{ij} \phi(e_j)$. If $\psi(\sum_{i=1}^d \alpha_i
e_i) = \sum_{i=1}^d \alpha_i e_i$, 
with $\alpha_i \in \pi^{-1} \aplus_F$,
then this translates into
$\psi(\sum_{i=1}^d \alpha_i q_{ij}) = \alpha_j$
for $1 \leq j \leq d$.

Let $\alpha_{i,n}$ be the coefficient of $\pi^n$ in $\alpha_i$, and
likewise for $q_{ij,n}$.
Since $\psi(1/\pi)=1/\pi$, the equations
$\psi(\sum_{i=1}^d \alpha_i q_{ij}) = \alpha_j$ then tell us that
for $1 \leq j \leq d$:
\[  \sum_{i=1}^d \alpha_{i,-1} q_{ij,0} = \phi(\alpha_{j,-1}). \] 
Since $N(V)/\pi N(V) \simeq \dcris(V)$ as filtered 
$\phi$-modules (cf the appendix to \cite{LBcr}), 
the above equations say that
$1$ is an eigenvalue of $\phi$ on $\dcris(V)$. It is easy to see
that if a representation has $\geq 0$ weights and $\dcris(V)^{\phi=1}
\neq 0$, then $V$ has a quotient isomorphic to $\Qp$. 
\end{proof}

\begin{rema}
It is proved in \cite{LBcr} that $\dcris(V) = (\bhol{,F} 
\otimes_F N(V))^{G_F}$ and that if $\on{Fil}^{-h} \dcris(V) =
\dcris(V)$, then $(t/\pi)^h \bhol{,F} 
\otimes_F \dcris(V) \subset \bhol{,F} 
\otimes_F N(V)$.

In all the above constructions, one could
therefore replace $\dhol(V)$ by $\bhol{,F} 
\otimes_F \pi^h N(V)$. For example, the image of $\Omega_{V,h}$
is included in $(\bhol{,F} 
\otimes_F \pi^h N(V))^{\psi=1}$.
\end{rema}

%%\Addresses
\end{document}